\documentclass{article}
\usepackage{amsthm}
\usepackage{amsmath}
\usepackage{amsfonts}
\usepackage{amssymb}
\usepackage{graphicx}
\usepackage{indentfirst}

\newtheorem{thm}{Theorem}[section]
\newtheorem{definition}[thm]{Definition}
\newtheorem{theorem}[thm]{Theorem}

\newtheorem{prop}[thm]{Proposition}
\newtheorem{lemma}[thm]{Lemma}
\newtheorem{conjecture}[thm]{Conjecture}
\newtheorem{corollary}[thm]{Corollary}
\begin{document}
 
\begin{center}{\bf \LARGE On a Basis for the Framed Link \\
\vspace{3pt} Vector Space Spanned by Chord \\ 
\vspace{3pt} Diagrams}
\bigskip

{\small
\begin{tabular}{lll} Bryan Bischof\mbox{\hspace{2pt}\footnotemark[1]} &Roman Kogan\footnotemark[1] & David N. Yetter\footnotemark[1]\\
Westminster College & SUNY Stony Brook & Kansas State University\\
New Wilmington, PA 15044 & Stony Brook, NY 11794 & Manhattan, KS 66506\\
bischobe@westminster.edu & romwell32@yahoo.com & dyetter@math.ksu.edu\end{tabular}}
\footnotetext[1]{The authors were partially supported by the Kansas State University REU and NSF grant GOMT530725}
\smallskip

14 January 2008\end{center}

\noindent{\bf Abstract} {\em In view of the result of Kontsevich, \cite{Kon} now often called ``the fundamental theorem of Vassiliev theory'', identifying the 
graded dual of the associated graded vector space to the space of Vassiliev invariants filtered by degree with the linear
span of chord diagrams modulo the ``4T-relation'' (and in the
unframed case, originally considered in \cite{Vas}, \cite{Kon},  and \cite{bn1}, the ``1T-'' or ``isolated chord relation''), it is a problem of some interest to provide a basis for the space of chord diagrams modulo the 4T-relation.

We construct the basis for the vector space spanned by chord diagrams with $n$ chords and $m$ link components, modulo $4$T relations for $n \leq 5$.}


\section{Introduction}

Perhaps due to Vassiliev's original formulation in terms of the space of knots \cite{Vas}, most work on Vassiliev theory has dealt with Vassiliev invariants of knots (e.g. \cite{bn1, Kon}), most often in the unframed setting.

In fact the entire development of the subject works equally well for links with any number of components, and (except for Kontsevich's transcendental methods that require the 1T-relation) in the framed setting. 

\begin{figure}[h]
\begin{center}
 \includegraphics[width = 12cm]{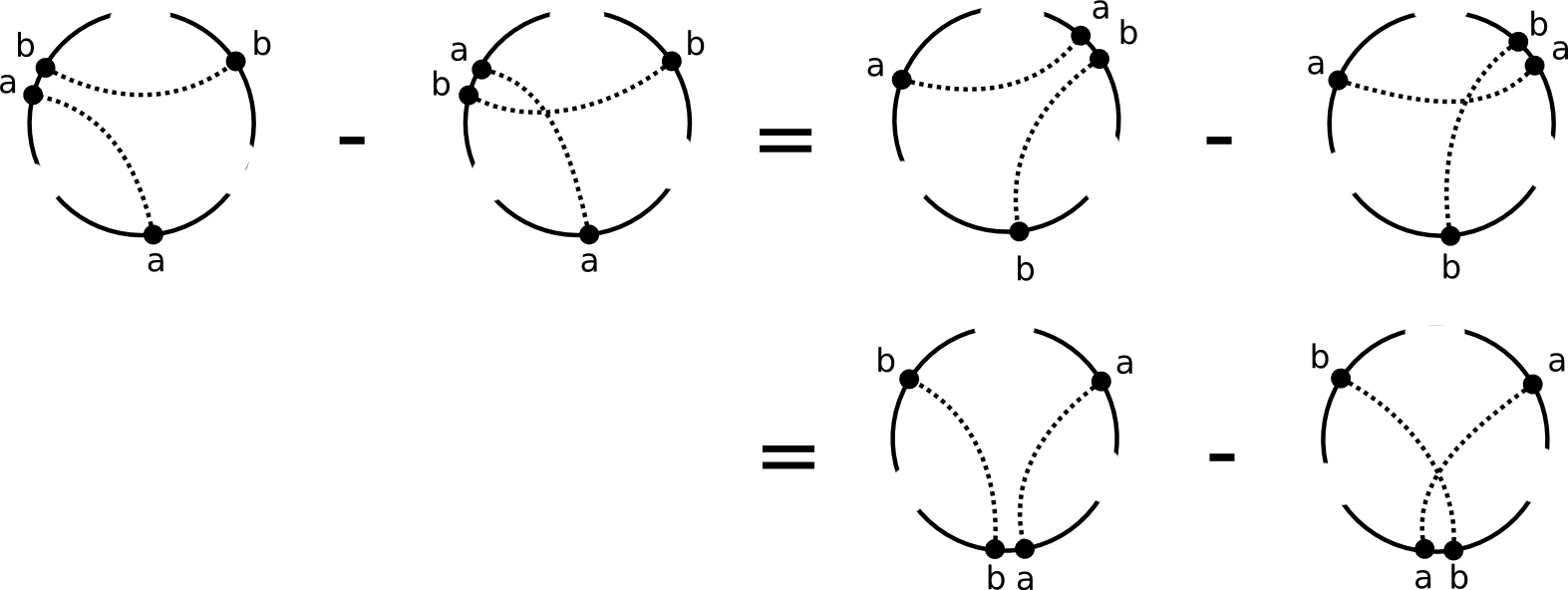}
 \caption{4T relations \label{4T}}
 \end{center}
\end{figure}

In particular, ${\cal V}^m$, the vector space of rational Vassiliev invariants of framed $m$-component links is naturally filtered by degree, and
the obvious generalization of Kontsevich's fundamental theorem
holds:

\begin{thm}

$\left({\cal V}_{n}^{m}/{\cal V} _{n-1}^{m}\right)$ is 
canonically dual to ${\cal A}_n^m$, the rational span of all
chord diagrams with $n$ chords and $m$ link components, modulo all instances of the 4T-relation. (See Figure \ref{4T}.)
\end{thm}

It is the purpose of this paper, in analogy to the work of Stanford \cite{Stan}, which gave an explicit basis for ${\cal A}_n^1$ for $n \leq 10$, to provide an explicit basis for
${\cal A}_n^m$ for $n \leq 5$ and all $m$.

We proceed as follows:  First we reduce the problem to the analogous problem for the subspace ${\cal C}_n^m$ spanned by chord diagrams whose underlying graph is connected, giving along the way a construction for the basis of ${\cal A}_n^m$ from
bases for ${\cal C}_k^m$, $k \leq n$. Then we describe a
data structure adequate for (a redundant) computer representation of all such chord diagrams, and give a canonical form theorem to identify a particular representation of each
distinct connected chord diagram.  We algorithmically generate all instances of 4T relations, and row reduce the resulting matrix. The generation of diagrams, comparison of diagrams, generation of instances of the 4T relation, and matrix reduction were all handled by computer. Specifically, programs written in Java 1.6.0 were used to the diagrams, converted to canonical form, and generate 4T relations. Mathematica was used to reduce the relations matrix and output the basis.

 Throughout we assume that the components of links are distinguishable from each other, and numbered $1, 2, \ldots, m$, though in the last section, we consider the question of passing to a basis for the space in which the components are indistinguishable by considering the
 action of the symmetric group ${\mathfrak S}_m$ on the $m$ components.
 
We do this because it is easy to see, even considering first order invariants, that there is a great deal more topological information in the Vassiliev invariants of framed links with distinguishable components:  if the components are indistinguishable, the space of first order Vassiliev invariants is spanned by the total framing number and the total linking number, while if the components are distinguishable, the entries of the linking matrix are first order Vassiliev invariants.
 
We use the standard convention that solid circles represents a component of the link, with chords drawn as dotted arcs, and when thinking of chord diagrams graph theoretically will simply refer to `solid arcs' and `dotted arcs'. We will refer to the place where a chord meets a solid arc as a `foot' of the chord. We will additionally assume the convention that all circles are considered to have clockwise orientation. As we most often will be considering framed knots and links, we will not in general consider the imposition of the 1T-relation. 

\section{Combinatorics} \label{combinatorics}

\begin{definition} The {\em degree} of a chord diagram is the number of chords.
\end{definition}

\begin{definition}
A chord diagram is {\em connected} if the underlying graph with solid circles as vertices and chords as edges is connected.
\end{definition}

Observe that the graphs here are what are usually called `generalized graphs':  loops and multiple edges are permitted.
 
\begin{definition}
A {\em full subdiagram} of a chord diagram is a chord diagram obtained by deleting some solid circles and those chords incident with them.
\end{definition}

\begin{definition}
A {\em connected component} of a chord diagram is a maximal connected full subdiagram.
\end{definition}

It is easy to see that every chord diagram is the disjoint union of its connected components.

We then have

\begin{prop} \label{4Tconn} For any four chord diagrams $D_1,\ldots, D_4$, related by an instance of the 4T relation, two solid circles lie in the same connected component of $D_1$ if and only if they lie in the same connected component of $D_i$ for all $i=1,\ldots, 4$. Moreover all but one of the connected components of the $D_i$'s are identical in all four diagrams, and the solid circles and the number of chords in the remaining component are the same in all four diagrams.
\end{prop}

\begin{proof}
``Only if'' is trivial, as is the second statement once the ``if'' is established.

For the other implication, observe that two solid circles $K$ and
$L$ lie in the same connected component if and only if there is a sequence 

\[ K = L_0, c_1, L_1, \ldots, c_n, L_n = L \]

\noindent in $D_1$ where each $L_j$ is a solid circle and each $c_j$ is a chord, and successive elements in the sequence are incident.  Without loss of generality, assume the sequence is minimal.

Now, if the instance of the 4T relation relating the $D_i$'s involves none of the $c_j$'s, the same sequence shows that $K$ and $L$ lie in the same component of each $D_i$.

Suppose instead that the instance of the 4T relation involved exactly one of the chords in the sequence, $c_j$.  In this case, one of the other $D_i$'s is obtained by exchanging the location of the place where $c_j$ meets one of the solid circles and the place where some chord adjacent $\kappa$ along the circle meets it, so again the same sequence of chords shows $K$ and $L$ lie in the same component. 

The other two of the $D_i$'s are obtained by moving the foot of $c_j$ adjacent to $\kappa$ to each side of the other foot of $\kappa$, which lies on some solid circle $\Lambda$.

In this case either

\[ L_0, c_1, L_1, \ldots, c_j, \Lambda, \kappa, L_j, \ldots, c_n, L_n = L\]

\noindent or 

\[ L_0, c_1, L_1, \ldots, c_{j-1}, \Lambda, \kappa, L_{j-1}, \ldots, c_n, L_n = L\]

\noindent shows that $K$ and $L$ lie in the same component, depending on whether $\kappa$ was incident with $L_j$ or $L_{j-1}$.

The remaining possibility is that the 4T relation involves two chords in the sequence.  By minimality, the $L_i$'s are distinct, and thus the two chords involved must be $c_j$ and 
$c_{j+1}$ for some $j$.  An argument like that above shows that
the same sequence of solid circles and chords verify that $K$ and $L$ lie in the same component of the diagram obtained by exchanging the adjacent feet of $c_j$ and $c_{j+1}$, while
the other in the other terms either

\[ L_0, c_1, L_1, \ldots, L_{j-1}, c_j, L_{j+1}, \ldots, c_n, L_n = L\]

\noindent or

\[ L_0, c_1, L_1, \ldots, L_{j-1}, c_{j+1}, L_{j+1}, \ldots, c_n, L_n = L\]

\noindent shows $K$ and $L$ are in the same component, depending on whether the foot of $c_j$ was moved near the other foot of $c_{j+1}$ or vice-versa. 
\end{proof}

Now, observe that the connected components of a chord diagram $D$ determine a partition $P(D)$ of the set of solid circles, which correspond to components of the link when the chord diagram arises as a summand in the Kontsevich integral, or as an element in the graded dual to the graded vector space associated to the filtered vector space of Vassiliev invariants of links on a fixed number of distinguishable components.

Any linear combination of chord diagrams $\sum_i a_i D_i$
can thus be expressed as a sum

\[ \sum_\Pi \sum_{P(D_i) = \Pi} a_i D_i \]

\noindent where the outer sum ranges over all partitions of the set of link components (or equiv. solid circles).  The terms can be further partitioned according to the number of chords in each connected component:

\[ \sum_\Pi \sum_{P(D_i) = \Pi} \sum_{\gamma(D_i) = \gamma}
a_i D_i \]

\noindent where the inners sum ranges over all compositions of the integer $n$ with as many summands as $\Pi$ has equivalence classes, and $\gamma(D_i)$ is the composition of $n$ whose $j^{th}$ summand is the number of chords in the $j^{th}$ connected component of $D_i$ when the connected components are
ordered by their lowest numbered solid circle.

Or, expressing the same observation in terms of the entire vector space ${\cal D}^m_n$ spanned by the set of chord diagrams on $m$ (distinguishable) solid circles, with $n$ (indistinguishable) chords, that

\[ {\cal D}^m_n = \bigoplus_\Pi {\cal D}^{m,\Pi}_{n,\gamma} \]

\noindent where ${\cal D}^{m,\Pi}_{n, \gamma}$ is spanned by all chord
diagrams on $m$ solid circles, with $n$ chords in which the partition of the set of solid circles determined by the connected components is $\Pi$, and the composition of $n$ determined by the number of chords in the components as in the
previous paragraph is $\gamma$.

It follows from Proposition \ref{4Tconn} that quotienting by 4T relations respects this direct sum decomposition. Thus we have

\begin{prop} \label{partitionsum}
If ${\cal A}^m_n$ is the quotient of ${\cal D}^m_n$ by all the subspace spanned by all instances of the 4T relation, it admits 
a corresponding decomposition 

\[ {\cal A}^m_n = \bigoplus_\Pi {\cal A}^{m,\Pi}_{n, \gamma} \]

\noindent where ${\cal A}^{m,\Pi}_{n, \gamma}$ is the quotient of ${\cal D}^{m,\Pi}_{n, \gamma}$ by all instances of the 4T relation.
\end{prop}

\begin{proof} It suffices to note that by Proposition \ref{4Tconn} all summands of a 4T relation lie in the same
${\cal D}^{m, \Pi}_{n, \gamma}$.
\end{proof}

Now, let ${\cal C}^m_n := {\cal A}^{m,{\bf 1}}_{n, n}$, where ${\bf 1}$ denotes the partition of the components with a single equivalence class, and the second $n$ represents the composition of $n$ as a single summand.  Thus ${\cal C}^m_n$ is the span of all connected chord diagram on the $m$ components with $n$ chords, modulo 4T relations. 


We then have another consequence of Proposition \ref{4Tconn}:

\begin{prop} \label{compositiontensor}

\[ {\cal A}^{m, \Pi}_{n, \gamma} \cong \bigotimes_{i=1}^{|\Pi|} 
{\cal C}^{|\Pi_i|}_{\gamma_i} \]

\noindent where $|\Pi_i|$ denotes the number of solid circles in $\Pi_i$, the $i^{th}$ connected component, when the components are ordered by their earliest numbered solid circle.
\end{prop}

\begin{proof}
Notice that the corresponding statement without quotienting by 4T relations
 
\[ {\cal D}^{m, \Pi}_{n, \gamma} \cong \bigotimes_{i=1}^{|\Pi|}
{\cal D}^{|\Pi_i|, {\bf 1}}_{\gamma_i, \gamma_i} \]
 
\noindent is just a 
linearization of the statement that any chord diagram is a disjoint union of its connected components.  Proposition \ref{4Tconn} then implies that any 4T relation with summands in 
${\cal D}^{m, \Pi}_{n, \gamma}$ corresponds to a 4T relation in exactly one of the tensorands, from which the desired result follows.
\end{proof}

Propositions \ref{partitionsum} and \ref{compositiontensor} reduce the proof of the following theorem to a few details.

\begin{theorem}

Let $A^m_n$ denote the dimension of ${\cal A}^m_n$, defined above, and let $C^m_n$ denote the dimension of ${\cal C}^m_n$. Then we have:

\begin{equation} \label{connectedtoall}
 A_{n}^m=\sum_{c=1}^{m}\dfrac{1}{c!}\left[\sum_{\substack{m_1+\dots+m_c=m \\ 1\leq m_i\leq m}}\binom{m}{m_1,\dots,m_c}\sum_{\substack{n_1+\dots+n_c=n \\ 0\leq n_i\leq n}}\prod_{i\leq c}C^{m_i}_{n_i}\right]
\end{equation}
 
 \noindent Note that
 $C^r_s=0$ for $s<r-1$ and $C^1_0=1$ .
\end{theorem}

\begin{proof}

Now, by Propositions \ref{partitionsum} and \ref{compositiontensor} we have

\[ {\cal A}^m_n = \bigoplus_{\Pi, \gamma} {\cal A}^{m,\Pi}_{n,\gamma} \cong \bigoplus_{\Pi, \gamma} \bigotimes_{i=1}^{|\Pi|} {\cal C}^{|\Pi_i|}_{\gamma_i} \]

The product for the $C^{m_i}_{n_i}$'s is plainly the 
dimension of the tensor product in the above expression when the
$m_i$'s and $n_i$'s are the corresponding parts of the partition of the solid circles, $\Pi$, and the composition of $n$, $\gamma$.

It thus remains only to see that the summation and coefficients actually correspond to the direct sum decomposition.

Obviously the dimensions satsify

\[ A^m_n = \sum_{\Pi, \gamma} {\cal A}^{m,\Pi}_{n,\gamma} = \sum_{\Pi, \gamma} \prod_{i=1}^{|\Pi|} C^{|\Pi_i|}_{\gamma_i} \]

\noindent where the sum ranges over all partitions of the set
of solid circles, and all compositions of the number of chords
into the same number of summands as the partition has equivalence classes.

This in turn can be rewritten as 

\[ A^m_n = \sum_{c=1}^n \sum_{\Pi, \gamma} \prod_{i=1}^{|\Pi|} C^{|\Pi_i|}_{\gamma_i} ,\]

\noindent where the inner sum now ranges over all partitions and compositions into $c$ parts (equivalence classes or summands).

It thus suffices to count the pairs of a partition of the solid circles and a composition of $n$ as above in which the sizes of the equivalence classes (resp. corresponding summand) are
$|\Pi_1|, \ldots, |\Pi_c|$ (resp. $\gamma_1, \ldots, \gamma_c$) when the equivalence classes are ordered by the lowest numbered solid circle.  

This is difficult to do directly, but observe that if we partition the set of solid circles into an {\em ordered} family of disjoint subsets of cardinalties $m_1, \ldots, m_c$, and use the {\em ordered} sum $n_1 + \ldots + n_c = n$ to assign numbers of chords to each, each of the pairs we wish to count will occur in the list exactly $c!$ times.

Noting that there are 

\[ \left(\begin{array}{c} m \\ m_1, \ldots, m_c \end{array} \right) \]

\noindent partitions of an $m$ element set into an ordered family of disjoint subsets of sizes $m_1, \ldots, m_c$ suffices to establish the formula.
\end{proof}

\section{Computer Representations and \\ Canonical Form}
We can represent chord diagrams on a knot using Gauss codes (cf. \textbf{\cite{cmu}}). For example $\left[012120\right]$ represent a chord diagram with three chords, labeled $0$, $1$, and $2$ so that beginning at some point on the knot, proceeding clockwise, the feet of the chords occuring the prescribed order.
 
In the case of chord diagrams with more than one solid circle, we must specify which chord feet lie on which component.  For instance $\left[0121|20\right]$ would denote a chord diagram with
two solid circles, three chords labeled $0$, $1$, and $2$, so that beginning at some point on the first solid circle traveling clockwise, the feet of chords 0, 1, 2, and 1 are encounted, in that order, while beginning at some point on the second solid circle, traveling clockwise the feet of chords 2 and 0 are encountered in that order.

More precisely:

\begin{definition}
A {\em string representation} of a chord diagram with $m$ solid circles and $n$ chords consists of an array of integers

\[ c_1, c_2, \ldots, c_{2n}, \] 

\noindent such that each integer $0, \ldots, n-1$ occurs
exactly twice among the $c_j$'s, together with a set of indices

\[ \{ 1=j_1 \leq j_2 \leq \ldots \leq j_m < j_{m+1} = 2n+1 \} \]

\end{definition}

The convention of defining $j_{m+1} = 2n+1$ makes
$c_{j_i}, \ldots, c_{j_{i+1}-1}$ always be the indices of the
chords incident with the $i^{th}$ solid circle, listed in clockwise order from some starting point.

For computer representation, the pair of an array and a set of indices is implemented as indicated.  For human reading, we indicate the elements of the set of indices (other than 1) by putting the symbol $|$ in front of $c_{j_i}$ for each $i = 2\ldots m$.

Now, observe that a given chord diagram admits many string
representations.  For instance $[0121|20]$, 
$[1020|21]$, $[1012|20]$, and $[0102|12]$ all represent the
same chord diagram.

Plainly permuting the labels on the chords and changing the
starting points from which the chord feet are listed on each
solid circle leaves the diagram unchanged, but changes the string representation.  Likewise, since every string representation is determined by an enumeration of the chords and a choice of starting point on each solid circle, any two string representations of the same chord diagram are related by a permuation of the chord labels and a cyclic permution of the chord feet on each component.

Of course, string representations (with a fixed set of indices) may be ordered lexicographically. We will refer to the lexicographically earliest string representation of fixed chord diagram as its {\em canonical form}.

Now, plainly if we fix the starting points, on each component, the lexicographically earliest string representation with those starting points is obtained by renumbering the chords in order of their first occurrence in the string.  Once this is done, any string representation of a chord diagram in which their are chords between solid circles can be moved earlier in lexicographic order by moving the starting point so that the
first chord foot listed is the earliest numbered chord from an earlier numbered solid circle.

Given a string representation, the canonical form of the corresponding chord diagram may thus be found algorithmically as follows:  

For each solid circle in order beginning with the first, modify the string reprsentation as follows:
 
If the solid circle is incident with chords from earlier solid circles, perform a cyclic permutation of the chord labels on the component to place the lowest numbered chord label first, then renumber all the chords from left to right according to the first occurence.  

Otherwise, if the solid circle is incident with no chord from an earlier solid circle, generate all string representations obtained by cyclic permuation of the chord labels on the circle, and renumbering all chords from left to right, and select the lexicographically earliest.

The canonical form provides a convenient way to represent chord diagrams on a computer, and for comparison of chord diagrams whether by computer or by hand.

Given an string representation $S$, we will denote the corresponding canonical form by $c(S)$.

Canonical form also provides a way of generating an exhaustive list of all chord diagrams with given numbers of solid circles, $m$, and chords, $n$:  G

Initialize an empty set of 'retained' string representations.

Generate in sequence all ordered sets of indices
$1 = j_1 \leq j_2 \leq \ldots \leq j_m < j_{m+1} = 2n+1$.  For each sequence of indices, generate (in lexicographic) order all sequences $c_1, \ldots, c_{2n}$ of integers in which each of the integers $0, \ldots, n-1$ occur exactly twice.  As each sequence is generated, compare the the canonical form $c(S)$ of the string representation $S = ([c_1, \ldots, c_{2n}], \{j_1, \ldots j_{m+1}\})$ with each of the retained string representations.  If it is not equal to any of them, retain $c(S)$ and continue, otherwise discard $c(S)$ and continue.

When all sets of indices and sequences have been generated, the set of retained string representations contains exactly the canonical form of each chord diagram one $m$ solid circles and $n$ chords.

From these, the set of all canonical forms of connected chord diagrams can be extracted by applying any standard spanning forest algorithm to the graph whose vertices are the solid circles and whose edges are the chords, and discarding any in which the tree generated with the first solid circle as root does not contain all $m$ solid circles.

Finally, given the set of canonical forms of connected chord diagrams on $m$ solid circles and $n$ chords, we can produce a complete, though redundant, list of all 4T relations among them as follows:

Given a canonical form 

\[ S = c(S) = ([c_1, \ldots, c_{2n}], \{j_1, \ldots j_{m+1}\}),  = [| c_1 \ldots | \ldots | \ldots c_{2n}|] \]
 
\noindent for each $i$ and each $k$
such that 

\[ j_i \leq k < k+1 < j_{i+1} \]

If $a = c_k \neq c_k+1 = b$ then there are two 4T relations

\begin{eqnarray*}
 [ \ldots a b \ldots a \ldots b \ldots ] -
c([ \ldots b a \ldots a \ldots b \ldots ]) - & & \\
c([ \ldots b \ldots b a \ldots b \ldots ]) +
c([ \ldots b \ldots a b \ldots b \ldots ]) & = & 0 \end{eqnarray*} 

\noindent and

\begin{eqnarray*}
 [ \ldots a b \ldots a \ldots b \ldots ] -
c([ \ldots b a \ldots a \ldots b \ldots ]) - & & \\
c([ \ldots a \ldots a \ldots b a \ldots ]) +
c([ \ldots a \ldots a \ldots a b \ldots ]) & = & 0 \end{eqnarray*}  

\noindent ({\em mutatis mutandis} when the other chord feet labeled $a$ and $b$ occur in different positions in the string representation relative to each other and the adjacent pair $a = c_k, b = c_{k+1}$ in the string representation).  Here the other symbols in the human-readable presentation of the string representation remain in the same order.  We leave it as an exercise to the reader to describe how the $j_k$'s change.

Likewise given a canonical for $S$ with $j_i < j_{i+1}-1$
and $b = c_{j_i} \neq c_{j_{i+1}-1} = a$ there are 
two 4T relations

\begin{eqnarray*}
 [ \ldots | b \ldots a | \ldots a \ldots b \ldots ] -
c([ \ldots | a \ldots b | \ldots a \ldots b \ldots ]) - & & \\
c([ \ldots | b \ldots | \ldots \ldots b a \ldots ]) +
c([ \ldots | b \ldots | \ldots \ldots a b \ldots ]) & = & 0 \end{eqnarray*}  

\noindent and

\begin{eqnarray*}
 [ \ldots | b \ldots a | \ldots a \ldots b \ldots ] -
c([ \ldots | a \ldots b | \ldots a \ldots b \ldots ]) - & & \\
c([ \ldots |  \ldots a | \ldots b a \ldots b \ldots ]) +
c([ \ldots |  \ldots a | \ldots a b \ldots b \ldots ]) & = & 0 \end{eqnarray*}

\pagebreak

\section{Connected Basis}

Having algorithmically determined an exhaustive list of the connected chord diagrams for a given number of solid circles $m$ and chords $n$, and an exhaustive (though redundant) list of 4T relations applying to them, we generate a matrix whose columns are indexed by the connected diagrams, and rows give the 4T relations.  Typically each row has four non-zero entries, two $1$'s and two $-1$'s, though sometimes coincidences of canonical forms cancel some entries, or give entries of $2$ or $-2$.

\begin{figure}[h]
\begin{center}
 \includegraphics[width = 3cm]{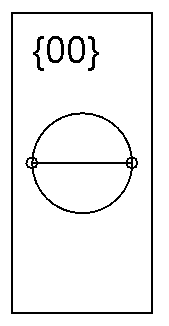}
 \caption{Basis of ${\cal C}_1^1$ \label{C_1^1}}
 \end{center}
\end{figure}

\begin{figure}[h]
\begin{center}
 \includegraphics[width = 6cm]{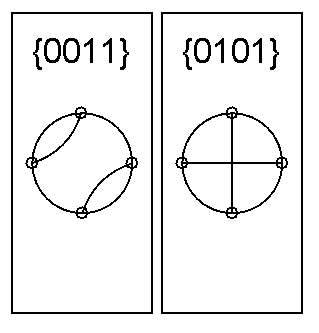}
 \caption{Basis of ${\cal C}_2^1$ \label{C_2^1}}
 \end{center}
\end{figure}

\begin{figure}[h]
\begin{center}
 \includegraphics[width = 9cm]{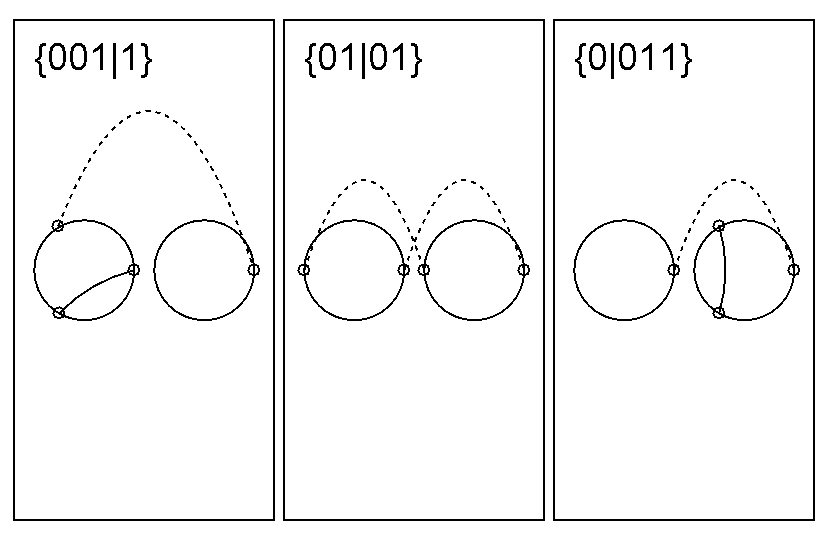}
 \caption{Basis of ${\cal C}_2^2$ \label{C_2^2}}
 \end{center}
\end{figure}

\begin{figure}[h]
\begin{center}
 \includegraphics[width = 9cm]{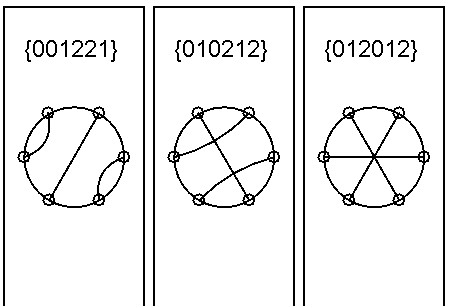}
 \caption{Basis of ${\cal C}_3^1$ \label{C_3^1}}
 \end{center}
\end{figure}

\begin{figure}[h]
\begin{center}
  \includegraphics[width = 4cm]{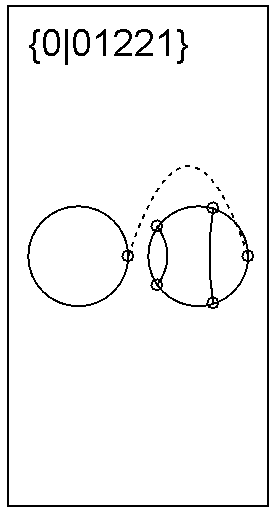}
 \caption{Basis of ${\cal C}_3^2$ (conclusion) \label{C_3^2}}
 \end{center}
\end{figure}

\begin{figure}[h]
\begin{center}
 \includegraphics[width = 12cm]{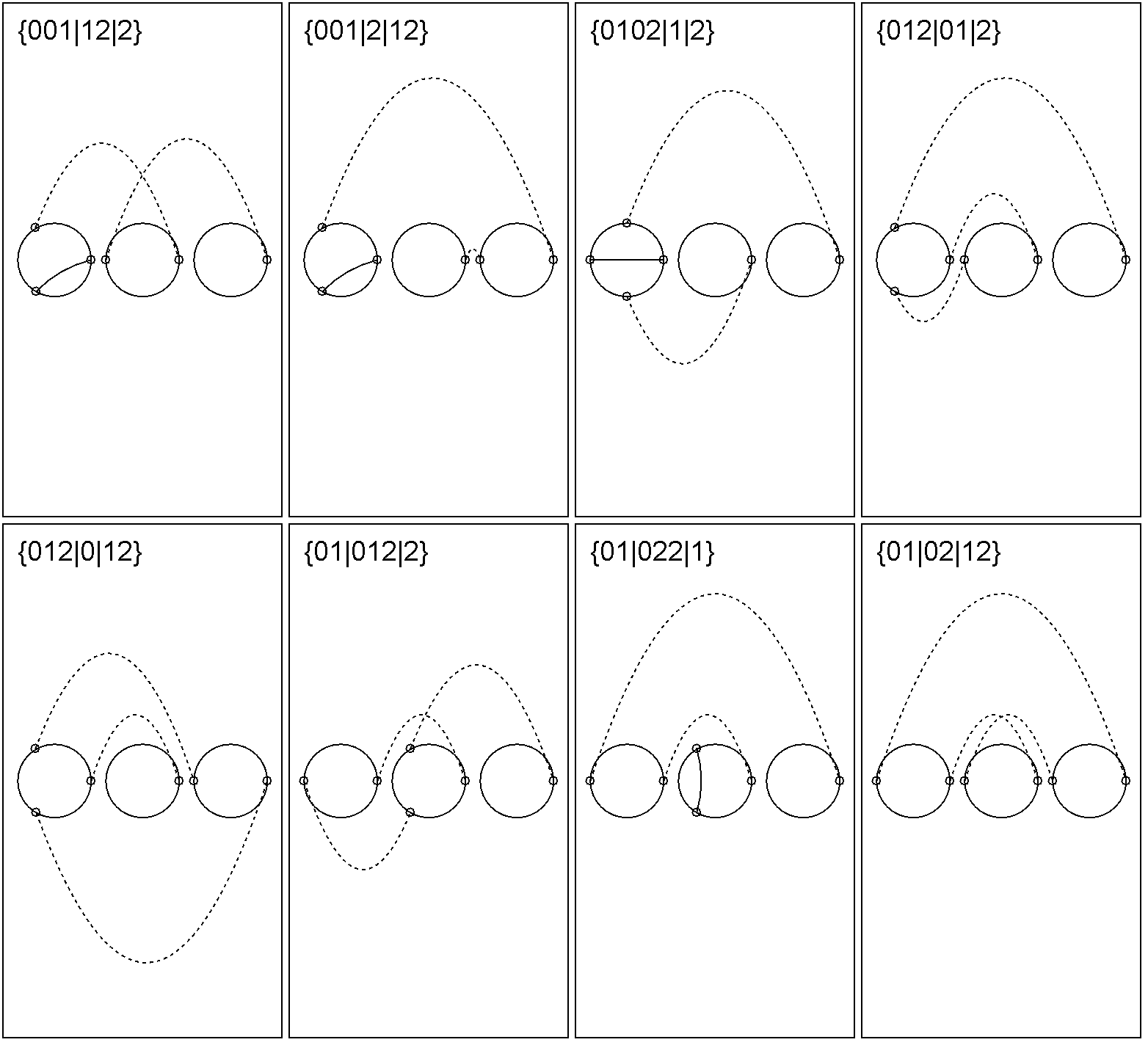}
 \caption{Basis of ${\cal C}_3^3$ (beginning) \label{C_3^3.0}}
 \end{center}
\end{figure}

\begin{figure}[h]
\begin{center}
  \includegraphics[width = 12cm]{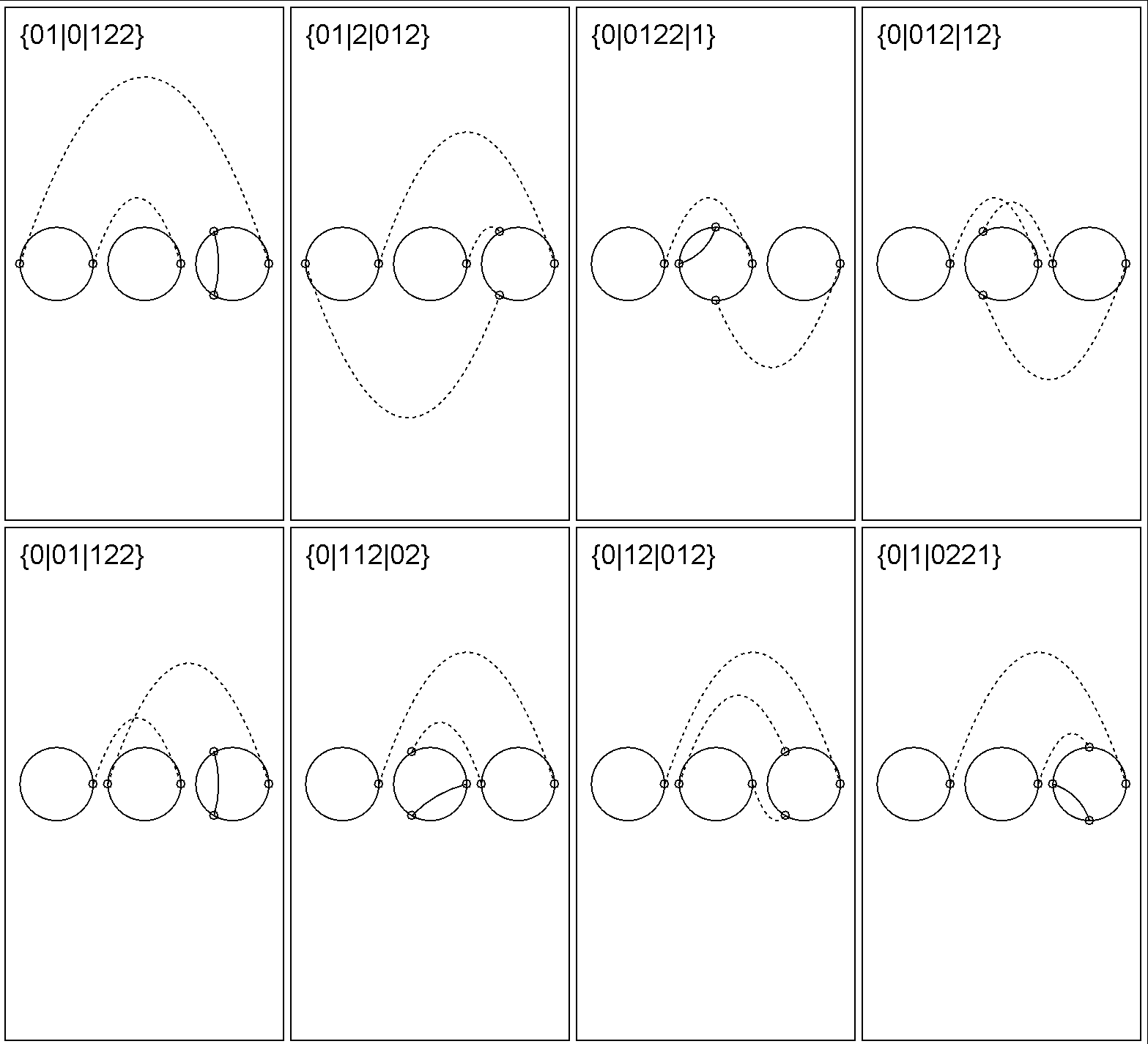}
 \caption{Basis of ${\cal C}_3^3$ (conclusion) \label{C_3^3.1}}
 \end{center}
\end{figure}

Row reducing the matrix then gives a set of linearly independent relation, the rows of the corresponding reduce row echelon form matrix. Now observe that the diagrams corresponding to the pivot columns are linear combinations of the other diagrams with non-zero coefficients in the same row, while all relations involving a diagram corresponding to a non-pivot column are either zero, or involve a pivot diagram. Hence, the set of diagrams minus those corresponding to the pivot columns will be a basis for the space ${\cal C}^m_n$.  Table \ref{C^m_n} gives the dimensions up to $n=5$.  In Figures \ref{C_1^1} to \ref{C_3^3.1} we show explicitly the chord diagrams in the basis computed for $n \leq 3$, and $m \leq n$.  A .zip archive of
figures in .png format containing the chord diagrams in the basis computed for $n \leq 5$, $m \leq n$, and $n=1$, $m=6$ can be downloaded from
http://www.math.ksu.edu/main/events/KSU-REU/BasisDiags.zip .  The case of ${\cal C}_n^{n+1}$ is discussed for all $n$ in Section \ref{equivariantbases} below.

\begin{table}[h]
\caption{$C^m_n$ \label{C^m_n}}
\centering
\begin{tabular}{c | c c c c c c | c}
 \hline\hline
${}_n{}^m$ & 1 & 2 & 3 & 4 & 5 & 6 & Total\\ [0.5ex]
\hline
1 & 1 & 1 &  &  &  &  & 2\\
2 & 2 & 3 & 3 &  &  &  & 8\\
3 & 3 & 9 & 16 & 16 &  & & 44\\
4 & 6 & 22 & 67 & 127 & 125 &  & 347\\
5 & 10 & 55 & 229 & 699 & 1347 & 1296 & 3636\\ 
\hline
\end{tabular}
\end{table}

\begin{figure}[h]
\begin{center}
 \includegraphics[width = 12cm]{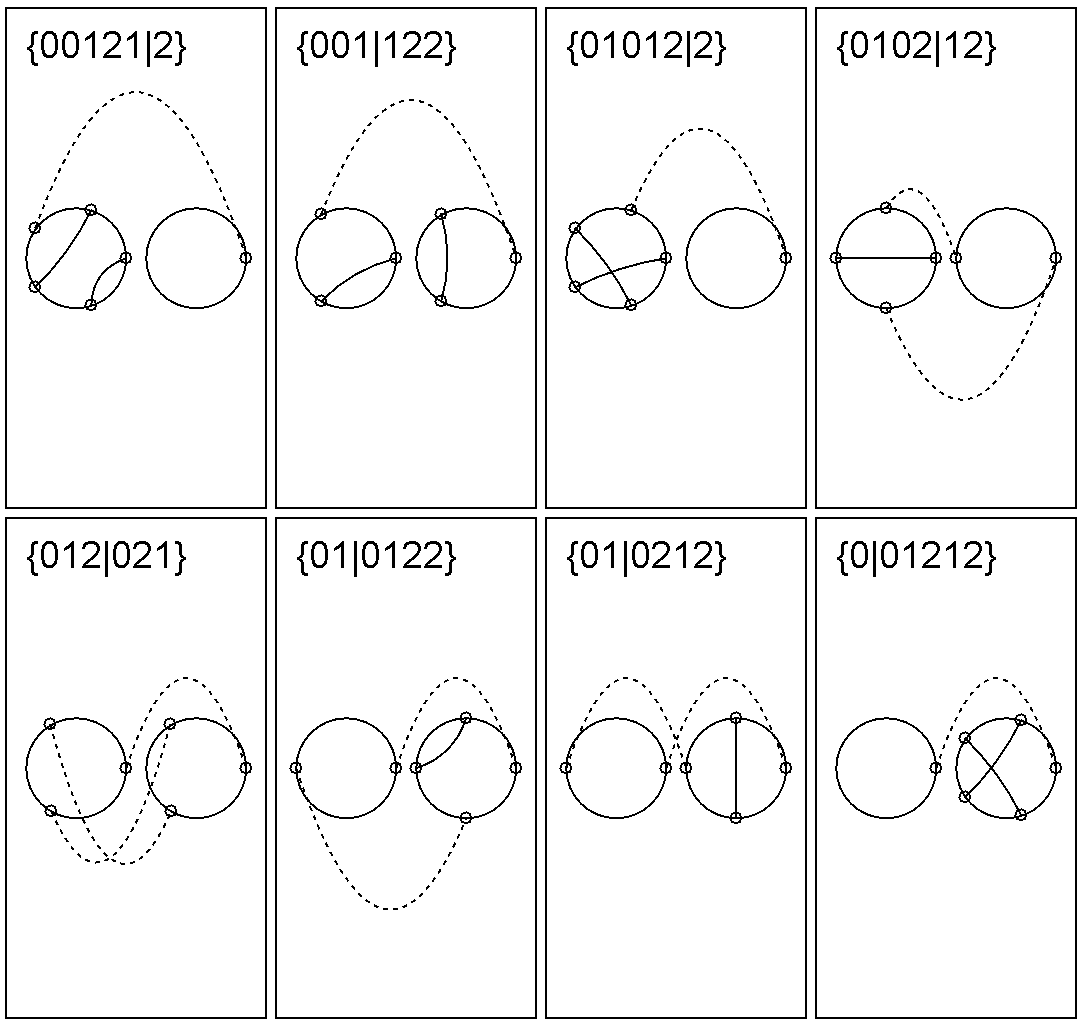}
\caption{Basis of ${\cal C}_3^2$ (beginning) }
 \end{center}
\end{figure}

\clearpage
\section{Complete Basis \label{completebasis}}
 With the values calculated for the connected basis we can use equation $\left(1\right)$ to generate the disconnected basis. Table \ref{A^m_n} give the dimensions $A_{n}^m$ for $m \leq 6$ and $n \leq 5$.\\
 
\begin{table}[h]
\caption{$A^m_n$ \label{A^m_n}}
\centering
\begin{tabular}{c | c c c c c c}
 \hline\hline
${}_n{}^m$ & 1 & 2 & 3 & 4 & 5 & 6\\ [0.5ex]
\hline
1 & 1 & 3 & 6 & 10 & 15 & 21 \\
2 & 2 & 8 & 24 & 59 & 125 & 237 \\
3 & 3 & 19 & 80 & 276 & 815 & 2088\\
4 & 6 & 44 & 241 & 1105 & 4340 & 14486 \\
5 & 10 & 99 & 682 & 3921 & 19468 & 81149 \\
\hline
\end{tabular}
\end{table}

It is easy to specialize Equation \ref{connectedtoall} to provide a general formula for $A_n^m$ for fixed values of $n$:

\begin{eqnarray*}
{\cal A}^m_1&=&\dfrac{m^2+m}{2}\\
\\
{\cal A}^m_2&=&\dfrac{m^4+3m^2}{8}+\dfrac{m^3+5m}{4}\\
\\
{\cal A}^m_3&=&\dfrac{m^6-287m^4}{144}+\dfrac{19m^5+325m^3}{48}-\dfrac{433m^2}{72}+\dfrac{23m}{6}\\
\\
{\cal A}^m_4&=&\dfrac{m^8-46375m^4}{384}+\dfrac{17m^7+26651m^3}{96}-\dfrac{209m^6}{64}+\dfrac{113m^5}{4}-{}\\&&\\&& {}-\dfrac{9775m^2}{32}+\dfrac{3107m}{24}\\
\\
{\cal A}^m_5&=&\dfrac{m^{10}+13188691m^5}{3840}+\dfrac{29m^9-151305m^6}{256}-\dfrac{1421m^8+23495m^7}{384}-{}\\&&\\&& {}-\dfrac{1139009m^4}{96}+\dfrac{4492697m^3}{192}-\dfrac{1897287m^2}{80}+\dfrac{557411m}{60}\\
\\
\end{eqnarray*}

\newpage

\section{Equivariant Bases \label{equivariantbases}}

Thus far we have considered the combinatorial underpinnings for the Vassiliev theory for framed links with distinguishable components.  One way of approaching the theory for links with indistiguishable components is to consider the action of the symmetric group on the set of link components that is induced on the space of Vassiliev invariants and on the space of chord
diagrams (modulo the 4T-relation).

In general the basis we computed in Section \ref{completebasis}
is not closed under the group action (for instance the basis
given above for ${\cal C}^2_3$ contains only $[01|0122]$, but not its image under the transposition of the solid circles, $[0012|12]$.  This is hardly surprising:  we started with a
spanning set on which ${\mathfrak S}_m$ acted, then reduced to
a basis.  In general in such a circumstance, the basis will not be a disjoint union of orbits--indeed often cannot be, no matter what the choice of spanning set (for example, given a non-trivial group character for a non-trivial group, the images of $1$ under the group action form a spanning set, but a basis for the underlying $1$-dimensional vector space will perforce contain only one of them).  When such a basis does exist, we call it an {\em equivariant basis}.


Nonetheless, as there is plainly a filtered vector space of Vassiliev invariants of links with indistinguishable components, and the same analysis will reduce the graded dual of its associated grade vector space to a vector space spanned by chord diagrams, we make

\begin{conjecture} \label{connectedconjecture}
 There exists an equivariant basis for ${\cal C}^m_n$ for all $m$ and $n$.
\end{conjecture}


From this the general conjecture

\begin{conjecture}
There exists and equivariant basis for ${\cal A}^m_n$ for all
$m$ and $n$ 
\end{conjecture}

\noindent follows from a combinatorial analysis similar to that in Section \ref{combinatorics}.

Even if we are interested in the case of distinguishable components, for ease of display and representation, it would be desirable to have an equivariant basis, as then the entire basis could be specified by giving a representative of each orbit.

Conjecture \ref{connectedconjecture} holds trivially for $m =1$,
and can be shown by {\em ad hoc} reasoning in the case $n = m = 3$:

\begin{prop}
${\cal C}^3_3$ admits an equivariant basis.
\end{prop}

\begin{proof}
The basis given in Figures \ref{C_3^3.0} and \ref{C_3^3.1} is already 
equivariant, being the union of two 6 element orbits (one in which the
underlying graphs are three vertex trees with a loop added to one leaf; the
other in which the underlying graph is a 2-cycle with an edge connecting
the third vertex to one of the verstices of the cycle); a three element
orbit (in which the underlying graph is a three vertex tree with a loop 
attached to the non-leaf); and a fixed point (with the 3-cycle as underlying
graph).
\end{proof}

In the next section, we prove the special case of Conjecture \ref{connectedconjecture} in which $m = n+1$, and with it give a description of the
basis for ${\cal C}^{n+1}_n$ for all $n$.

In Section \ref{C^2_n} we prove the conjecture holds fpr $m = 2$,
proving along the way a lemma that would appear to be useful for the proof in general.
\subsection{${\cal C}^{n+1}_n$}
Connected chord diagrams in ${\cal C}^{n+1}_n$ are special in that their underlying graphs are tree graphs.  We then have

\begin{theorem}
If $D_1$ and $D_2$ are two connected chord diagrams in ${\cal D}^{n+1}_n$
and $D_2$ is obtained from $D_1$ by permuting the feet of  chords incident with a given solid circle, then the images of $D_1$ and $D_2$ in ${\cal C}^{n+1}_n$ are equal.
\end{theorem}

For example, in the case of transposing two chords:

\begin{center}
 \includegraphics[height = 4cm]{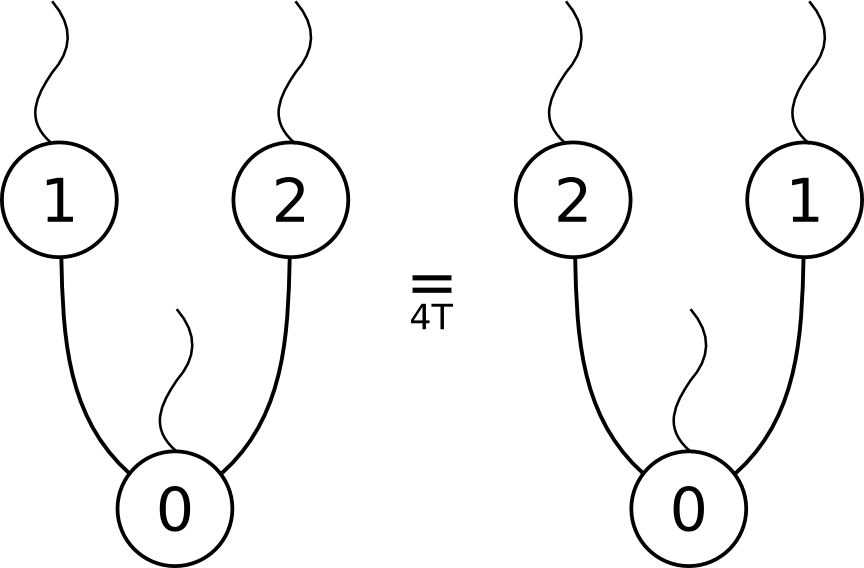}
\end{center}

\noindent where the wavy lines indicate an arbitrary number
of chords.

\begin{proof}
As any permutation can be written as a product of transpositions
of adjacent elements, it plainly suffices to show
 
\begin{lemma}

If $D_1$ and $D_2$ are two connected chord diagrams in ${\cal D}^{n+1}_n$
and $D_2$ is obtained from $D_1$ by transposing the adjacent
feet of a
pair chords incident with a given solid circle, then the images of $D_1$ and $D_2$ in ${\cal C}^{n+1}_n$ are equal.
\end{lemma}

\begin{proof}
Make the underlying tree of the chord diagram into a rooted 
tree by chosing the solid circle, say $0$, on which the adjacent pair of chords were interchanged in passing from $D_1$ to $D_2$, to be the root.  Now, let $m$ be
the total number of descendants of the nodes $1$ and $2$ whose
edges (chords) by which they were adjacent to $0$ (or, equivalently the total number of chords in the rooted subtrees at $1$ and $2$.  Pictorially:

\begin{center}
 \includegraphics[height = 4cm]{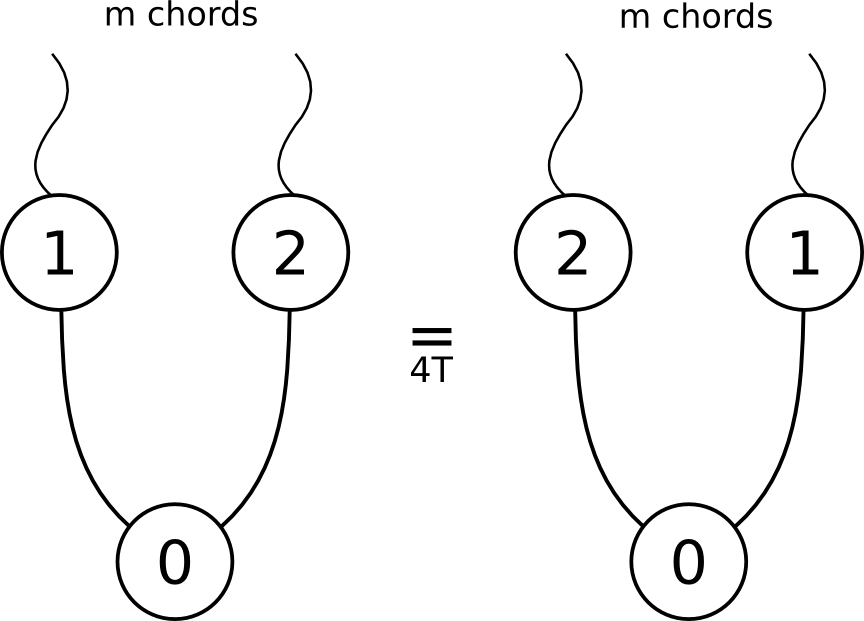}
\end{center}

We proceed by induction on $m$. For the base case let $m$ equal zero:

 \begin{center}
 \includegraphics[height = 4cm]{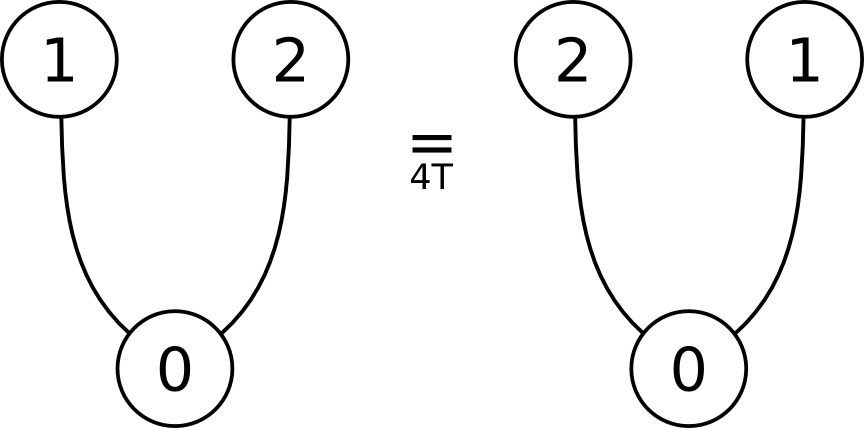}
\end{center}

This follows immediately from the 4T relation:

\begin{center}
 \includegraphics[width = 12cm]{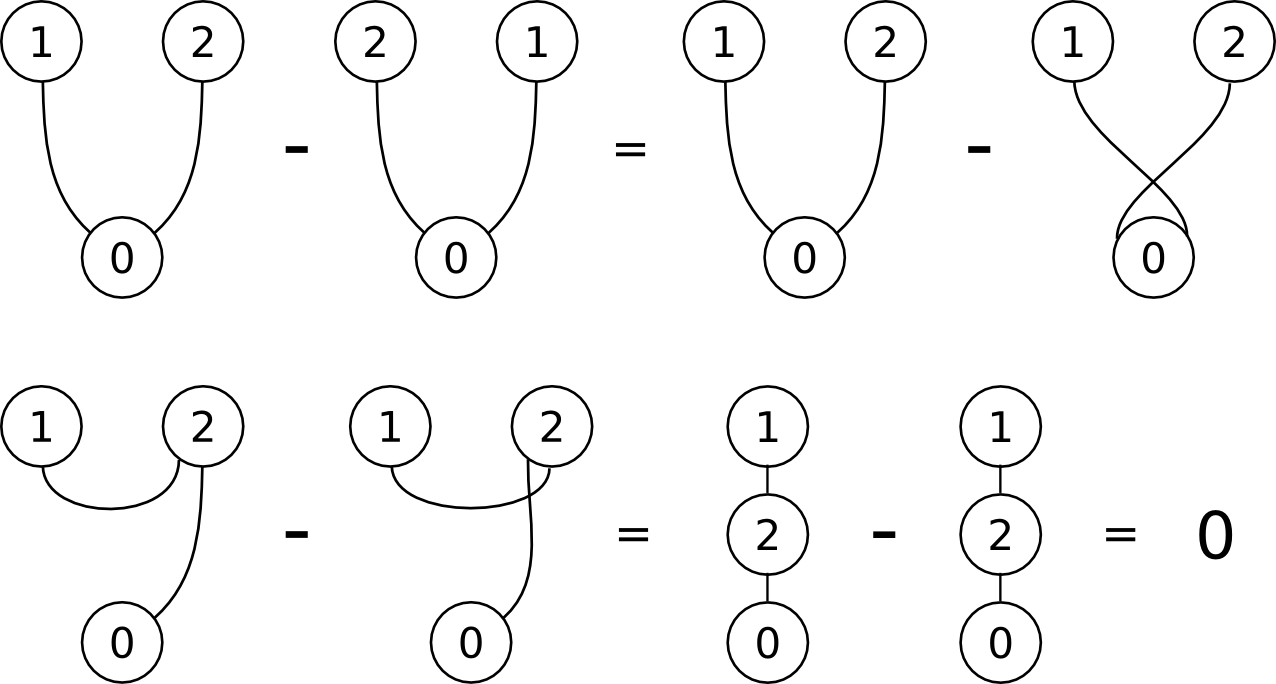}
\end{center}

Our inductive hypothesis is then that the lemma holds whenever $m < k$:

\begin{center}
 \includegraphics[height = 4cm]{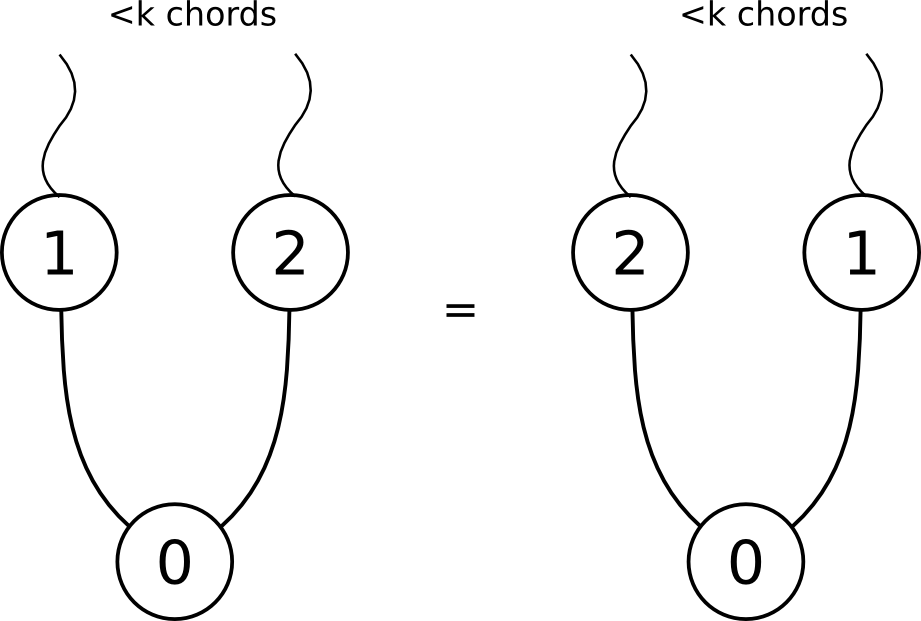}
\end{center}

Now consider a pair of adjacent chords for which $m=k$:

 \begin{center}
 \includegraphics[height = 4cm]{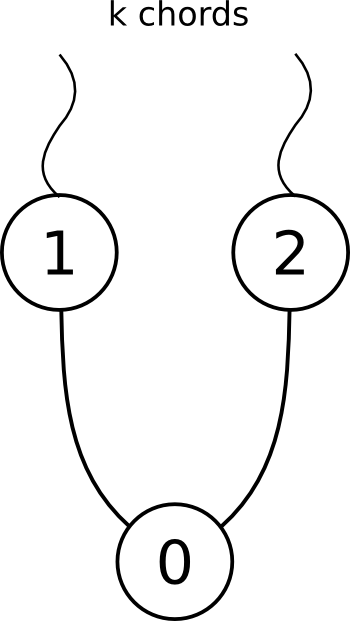}
\end{center}

Using our inductive hypothesis we proceed as follows:
\begin{center}
 \includegraphics[width = 12cm]{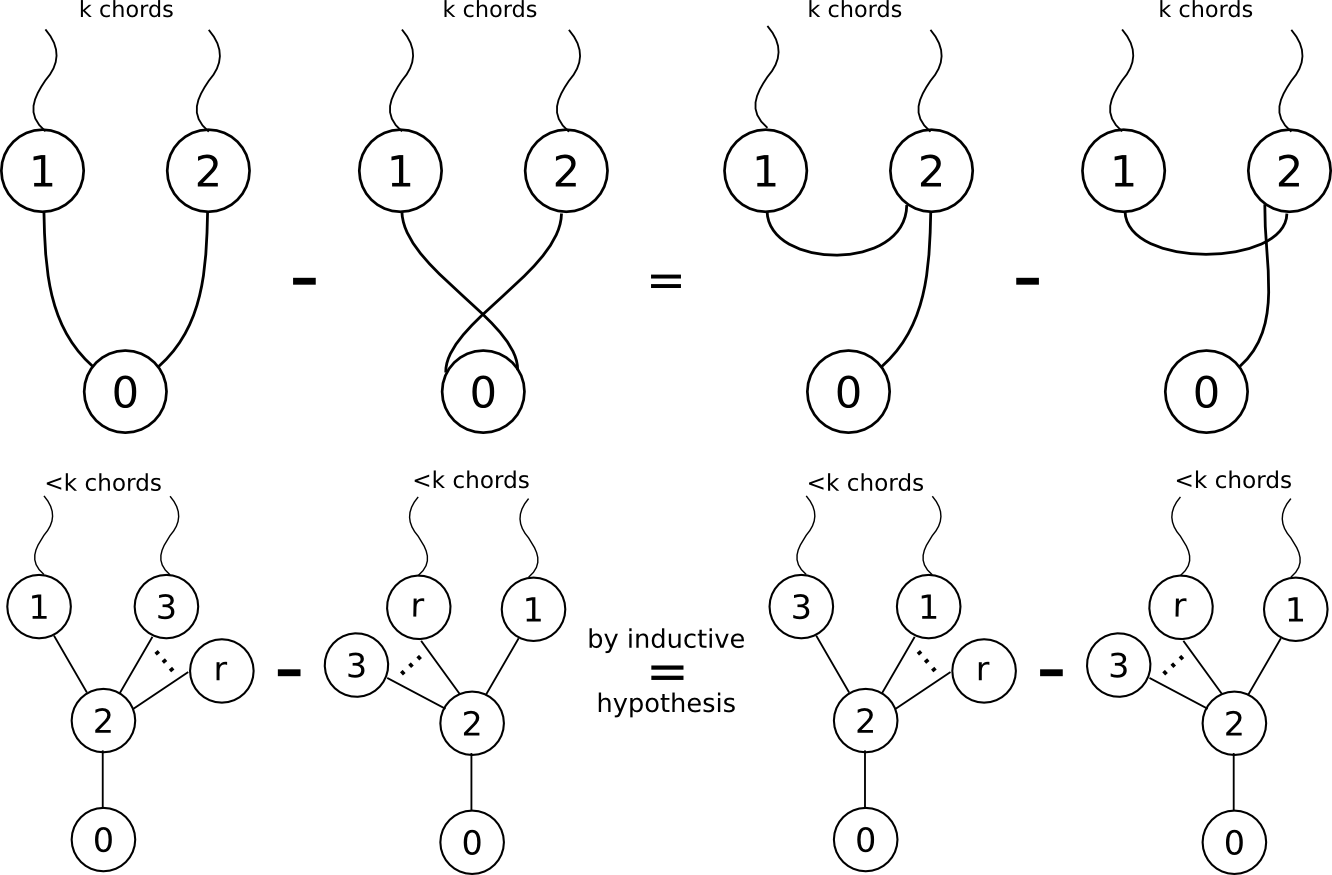}
\end{center}

\noindent Repeated application of the inductive
hypothesis with the root at $2$ then gives

\begin{center}
 \includegraphics[width = 12cm]{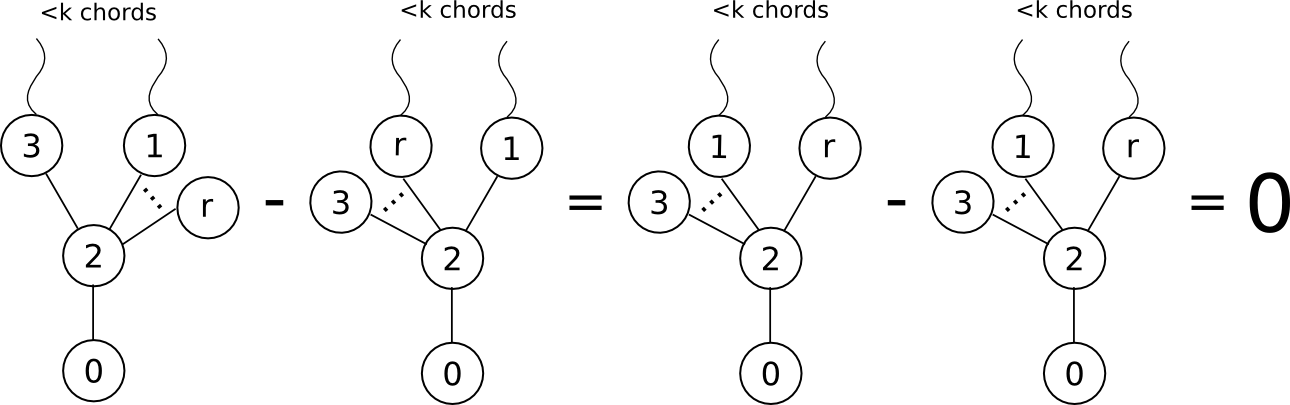}
\end{center}

\end{proof}

And thus, the theorem holds.
\end{proof}

We have then have:

\begin{corollary}
Any two connected chord diagrams $D_1$ and $D_2$ have the same labeled tree as underlying graph then $D_1$ and $D_2$ congruent modulo 4T relations.
\end{corollary}

Now observe that any instance of the 4T relation applied to a component with a tree as underlying graph, written as
in Figure \ref{4T}, as an equation between differences of
diagrams which differ only by swapping adjacent chord feet, asserts that two differences which are each zero by be previous Corollary are equal.
 
From this it follows that

\begin{theorem}The equivalence relation on $D^{n+1}_n$ induced by all instances of the 4T relation is the same as the equivalence relation given by linearizing the equivalence relation on chord diagrams in which two chord diagrams with the same underlying labeled graph (tree) are equaivalent.
\end{theorem}

Thus

\begin{theorem} \label{treebasis}
 Any basis of ${\cal C}^{n+1}_n$ consists of a set of chord diagrams containing exactly one representative for each isomorphism class of underlying labeled tree.
\end{theorem}

\noindent which result also completes our explicit construction of bases for ${\cal C}^m_n$ and ${\cal A}^m_n$ for $n \leq 5$.

It also follows that

\begin{theorem}
${\cal C}^{n+1}_n$ admits an equivariant basis.
\end{theorem}

\begin{proof}
 By the previous argument we have shown that every chord diagram in ${\cal D}^{n+1}_n$ is equal in ${\cal C}^{n+1}_n$ to a chord diagram in the basis, thus the basis given by Theorem \ref{treebasis} is the union of all orbits of images of chord
 diagrams in ${\cal C}^{n+1}_n$
 under the action of the symmetric group and thus {\em a fortiori} equivariant.
\end{proof}

\begin{corollary}
$C^{n+1}_n = \left(n+1\right)^{n-1}$.
\end{corollary}
\begin{proof}
 The number of labelled trees on $m$ vertices is $\left(m\right)^{m-2}$ by the Cayley-Borchardt Formula \cite{cla}. The number of vertices is the number of components, $n+1$.
\end{proof}

\subsection{Constructing Equivariant Bases, and the case of ${\cal C}^2_n$} \label{C^2_n}

In general, given a basis ${\cal B}^m_n$ of ${\cal C}^m_n$, the 
orbits of a basis vector will 
either be contained in the basis, or intersect the basis in a 
proper non-empty subset of the orbit.  In the former case, we will call
the orbit {\em complete}, in the later, {\em incomplete}.  If the orbit of any
basis vector is incomplete, we will say ${\cal B}^m_n$ {\em has incomplete 
orbits}.  We will denote the orbit of a vector $v$ by ${\mathfrak S}_m(v)$.

It is completely trivial to observe

\begin{prop}
A basis  ${\cal B}^m_n$ of ${\cal C}^m_n$ is equivariant if it does not
have incomplete orbits.
\end{prop}

To prove Conjecture \ref{connectedconjecture} either in general or for a particular $n$ and $m$ one need a construction which will reduce the number of
incomplete orbits.

The following lemma is useful:

\begin{lemma}
If $b \in {\cal B}^m_n$ lies in an incomplete orbit, and 
$\sigma(b) \not{\in} {\cal B}^m_n$ for some $\sigma \in {\mathfrak S}_m$, 
then, when we write $\sigma(b)$ as a linear combination of the basis elements

\[ \sigma(b) = \sum_{\beta \in {\cal B}^m_n} c_\beta \beta \]

\noindent  either $c_\beta \neq 0$ for some $\beta \in {\mathfrak S}_m(b) \cap {\cal B}^m_n$ or $c_\beta \neq 0$ for some $\beta$ in a different incomplete orbit.
\end{lemma}

\begin{proof}
Plainly the only other possibility is that the non-zero coefficients 
are all for basis vectors in complete orbits.  However in this case,
applying $\sigma^{-1}$ to both sides of the equation gives an expression
for $b$ as a linear combination of other basis vectors, contradicting linear 
independence.
\end{proof} 

We will call an incomplete orbit {\em type I} if it contains a translate of a
basis element satisfying the first 
condition, and {\em type II} if it contains a translate of a basis element
satisfying the second condition.  
(Notice it is possible for an incomplete orbit to be of both type I 
and type II.)  

It is now easy to show

\begin{theorem}
${\cal C}^2_n$ admits an equivariant basis for all $n$
\end{theorem}

\begin{proof}
Given a basis ${\cal B}^2_m$, with no incomplete orbits, we are done.  It 
therefore suffices to show that given any basis ${\cal B}^2_m$ it is possible
to replace it with a basis with fewer incomplete orbits.  This follows
directly from the previous lemma and the following two lemmas: 
\end{proof}

\begin{lemma}
Given an incomplete orbit of type I for a basis $B^2_n$, intersecting
the basis in
the singleton $\{b\}$.  Replacing $b$ with the average (or sum) of the 
elements in ${\mathfrak S}_2(b)$ gives a basis  ${\cal B}^{2 \prime}_n$
with one fewer incomplete orbit.
\end{lemma}

\begin{proof}
Trivial.  (Observe that the average or sum of elements in an orbit is
a fixed point of the action.)
\end{proof}

\noindent and

\begin{lemma}
Given an incomplete orbit of type II for a basis $B^2_n$, intersecting
the basis in the singleton $\{b\}$, if 

\[ \sigma(b) = c_\beta \beta + \mbox{\rm other basis elements} \]

\noindent where $\beta$ lies in a different incomplete orbit and $c_\beta 
\neq 0$,  then replacing $\beta$ with $\sigma(b)$ gives a basis
 ${\cal B}^{2 \prime}_n$ with two fewer incomplete orbits.
\end{lemma}

\begin{proof}
Trivial.
\end{proof}

\section{Prospects}

The results contained herein should be seen as providing a test-bed for conjectures about Vassiliev invariants of
(framed) links by allowing explicit calculations at low degree.  

It would appear to be feasible with increased computing power to use the naive methods of the present paper to find bases for
${\cal C}^6_m$, but for $n > 6$ a more sophisticated approach would appear to be necessary to render the computations feasible.

\pagebreak

\end{document}